\documentclass[11pt]{article}

\usepackage{amsmath,amssymb,amsfonts,amsthm}

\newcommand{\Hom}{\mbox{Hom}\,}
\newcommand{\Ext}{\mbox{Ext}\,}

\newcommand{\Spec}{\mbox{Spec}\,}

\newcommand{\Supp}{\mbox{Supp}\,}
\newcommand{\gr}{\mbox{grade}\,}
\newcommand{\depth}{\mbox{depth}\,}
\renewcommand{\dim}{\mbox{dim}\,}
\newcommand{\cd}{\mbox{cd}\,}

\newcommand{\pd}{\mbox{proj.dim}\,}
\newcommand{\id}{\mbox{inj.dim}\,}

\newcommand{\gd}{\mbox{Gdim}\,}

\newcommand{\uhom}{{\mathbf R}\Hom}
\newcommand{\utp}{\otimes^{\mathbf L}}
\newcommand{\ugamma}{{\mathbf R}\Gamma}
\renewcommand{\H}{\mbox{H}}

\newcommand{\T}{\mathrm}

\newcommand{\lo}{\longrightarrow}

\newcommand{\fa}{\frak{a}}

\newcommand{\fm}{\frak{m}}
\newcommand{\fp}{\frak{p}}

\begin{document}
\date{}
\title{\bf Generalized local cohomology
and the Intersection Theorem\footnotetext{2000 {\it Mathematics
subject classification.} 13D45, 13D22, 13D25, 13D05.}
\footnotetext{{\it Key words and phrases.} Intersection Theorem,
cohomological dimension, complexes of modules.}\footnotetext{The
research of the first author was supported by a grant from IPM
(No. 82130014).}\footnotetext{The second author was supported by a
grant from University of Tehran (No.
511/3/702).}\footnotetext{Email addresses: dibaeimt@ipm.ir and
yassemi@ipm.ir}}

\author{Mohammad T. Dibaei $^{\it ab}$ and Siamak Yassemi$^{\it
ca}$\\
{\small\it $(a)$ Institute for Studies in Theoretical Physics and
Mathematics, P.O. Box 19395--5746, Tehran, Iran}\\
{\small\it $(b)$ Department of Mathematics, Teacher Training
University, Tehran, Iran}\\
{\small\it $(c)$ Department of Mathematics, University of Tehran,
Tehran, Iran}} \maketitle

\begin{abstract}

\noindent Let $R$ be commutative Noetherian ring and let $\fa$ be
an ideal of $R$. For complexes $X$ and $Y$ of $R$--modules we
investigate the invariant $\inf{\mathbf R}\Gamma_{\fa}({\mathbf
R}\Hom_R(X,Y))$ in certain cases. It is shown that, for bounded
complexes $X$ and $Y$ with finite homology, $\dim Y\le\dim{\mathbf
R}\Hom_R(X,Y)\le\pd X+\dim(X\otimes^{\mathbf L}_RY)+\sup X$ which
strengthen the Intersection Theorem. Here $\inf X$ and $\sup X$
denote the homological infimum, and supremum of the complex $X$,
respectively.

\end{abstract}

\baselineskip=18pt

\vspace{.3in}

\noindent{\bf 0. Introduction}

\vspace{.2in}

Cohomological dimension of a module $M$ over a finite dimensional
Noetherian ring $R$ with respect to an ideal $\fa$ is denoted by
$\cd(\fa, M)=\sup\{i\in\Bbb Z|\H^i_{\fa}(M)\neq 0\}$ and has been
studied by Faltings [{\bf Fa}] and Huneke--Lyubeznik [{\bf HL}],
where $\H^i_{\fa}(M)$ is the $i$-th right derived functor of the
section functor $\Gamma_{\fa}(M)$. The authors, in [{\bf DY}],
introduced the notion of cohomological dimension, $\cd(\fa, X)$,
of a bounded to the left complex $X$ with respect to an ideal
$\fa$ to be $\sup\{\cd(\fa,\H_{\ell}(X))-\ell|\ell\in\Bbb Z\}$. It
is shown that $-\inf{\mathbf R}\Gamma_{\fa}(X)\le\cd(\fa, X)$ and
equality holds whenever $X$ has finite homology modules.

Let $M$ be a finite $R$--module with finite projective dimension
and let $N$ be a finite $R$--module. A consequence of the New
Intersection Theorem from [{\bf PS}] is the dimension inequality:
$\dim N\le\pd M+\dim (M\otimes_R N)$, cf. [{\bf R}; 8.4.4].
Actually the following inequalities hold (see [{\bf Y1}; Remark
3.8]):

\noindent {\bf (T1)} $\dim N\le\dim{\mathbf R}\Hom_R(M,N)\le\dim
M\otimes_R N+\pd M$.

In particular, if we replace $N$ by $R$ then we have

\noindent {\bf (T2)} $\dim R\le\dim{\mathbf R}\Hom_R(M,R)\le\dim
M+\pd M$.

Since the notions ``dimension'' and ``cohomological dimension with
respect to an ideal'' are closed in some sense, it is natural to
ask about the following question: ``Do we have the inequality
$\cd(\fa, {\mathbf R}\Hom_R(M,N))\le\cd(\fa, M\otimes_R N)+\pd M$
for finite $R$--modules $M$ and $N$ with $\pd M<\infty$?'' In this
paper we give a positive answer to this question.

In 1967 Auslander [{\bf A}] introduced a new invariant for finite
modules, a relative homological dimension, called the Gorenstein
dimension and denoted by $\gd(-)$. It is well--known that for all
finite $R$--modules $M$ there is an inequality $\gd M\le\pd M$ and
equality holds if $\pd M<\infty$. The Gorenstein dimension is
extended for complexes in [{\bf Y2}] (see also [{\bf C}]).

After introductory section 1, we consider complexes of finite
homologies and show in section 2 that if $X$ and $Y$ are bounded
complexes, then
$$\dim{\mathbf R}\Hom_R(X,Y)\le \dim R+\sup X-\inf Y$$provided
$\pd X$ or $\id Y$ is finite, or, $\gd X$ and $\pd Y$ are finite.
This result, which is a generalization of [{\bf HZ}; Theorem 3.2],
implies the equality $\dim{\mathbf R}\Hom_R(M,R)=\dim R$, provided
$\pd M<\infty$, which forces the left inequality of (T2) to be an
equality. This result motivates us to give some characterization
properties for a local ring $(R,\fm,k)$. For example, $(R,\fm,k)$
is Gorenstein (resp. Cohen--Macaulay) ring if and only if
$\dim{\mathbf R}\Hom_R(X,R)$ is finite (resp. $\dim{\mathbf
R}\Hom_R(X,L)$ is finite for some non--zero finite $R$--module
$L$) for all bounded complexes $X$ with finite homologies (see 2.3
and 2.4).

We next examine the Intersection Theorem and show that
$$\dim Y+\inf X\le\dim{\mathbf R}\Hom_R(X,Y)\le\dim(X\otimes^{\mathbf L}_RY)+\pd X+\sup X$$ for all bounded
complexes $X,Y$ with finite homologies and $\pd X<\infty$ (see
Theorem 2.6). This result generalizes (T1), which in particular
implies
$$\dim Y\le\dim{\mathbf R}\Hom_R(M,Y)\le\dim
M\otimes^{\mathbf L}_RY+\pd M$$ for all bounded complexes $Y$ and
all finite $R$-module $M$ with finite projective dimension (see
Corollary 2.7).

In section 3, we mainly investigate the cohomological dimension of
${\mathbf R}\Hom_R(X,Y)$ with respect to an ideal $\fa$ of $R$,
where $X,Y$ are appropriate $R$--complexes with $\pd X<\infty$.
Our main result is as follows:

If $X,Y$ are bounded complexes with finite homologies then

\[ \begin{array}{rl}
-\inf{\mathbf R}\Gamma_{\fa}(Y)+\inf X &\, \le -\inf{\mathbf
R}\Gamma_{\fa}({\mathbf R}\Hom_R(X,Y))\\
&\, \le -\inf{\mathbf R}\Gamma_{\fa}(X\otimes^{\mathbf L}_R Y)+\pd
X+\sup X,
\end{array} \]
where the first inequality holds under the extra condition $\Supp
Y\subseteq\Supp X$ (see Theorems 3.2 and 3.3).

Our results in section 3 show the connection between cohomological
dimension and the Intersection Theorem.

Throughout $R$ is a commutative Noetherian ring and $\fa$ is an
ideal of $R$.

\vspace{.3in}

\noindent{\bf 1. Notations}

\vspace{.2in}

An $R$-complex $X$ is a sequence of $R$-modules $X_\ell$ and
$R$-linear maps $\partial_\ell^X, \ell\in \mathbb{Z}$,
$$X=\cdots \lo X_{\ell+1} \stackrel{\partial^X_{\ell+1}}{\lo} X_\ell
\stackrel{\partial_\ell^X}{\lo} X_{\ell-1} \lo \cdots. $$ The
module $X_\ell$ is called the module in degree $\ell$, and the map
$\partial_\ell^X:X_\ell\lo X_{\ell-1}$ is the $\ell$-th
differential, and $\partial_\ell^X \partial_{\ell+1}^X=0$ for all
$\ell\in \mathbb{Z}$. An $R$--module $M$ is thought of as the
complex $M=0\lo M\lo 0.$

The {\em supremum} and  {\em infimum} of $X$ are defined by
\[ \begin{array}{rl} \sup \: X & =
\sup \: \{ {\ell} \in \Bbb Z | \H_{ \ell } (X) \neq 0 \}
\\ [.1in] \inf \: X & = \inf \: \{ {\ell} \in \Bbb Z | \H_{
\ell } (X) \neq 0 \},
\end{array} \]
denote $\sup X=-\infty$ and $\inf X=\infty$ if $\H_{\ell}(X)=0$
for all $\ell$; in this case $X$ is called homologically trivial.

A morphism $\alpha: X\lo Y$ is said to be a quasi-isomorphism, and
denoted by $X\simeq Y$, if the induced morphism
$\T{H}(\alpha):\T{H}(X)\lo \T{H}(Y)$ is an isomorphism. Support of
$X$ is defined by $\Supp X=\{\fp\in\Spec R|X_{\fp}\,\, \mbox{is
not homologically trivial}\}$.

The {\it derived category} ${\cal D}(R)$ is the category of
$R$--complexes localized at the class of all quasi--isomorphisms.
The full subcategories ${\cal D}_{+} (R)$, ${\cal D}_{-} (R)$, and
${\cal D}_{b} (R)$ of ${\cal D}(R)$ consist of complexes $X$ with
$\H_{\ell}(X)=0$ for, respectively, $\ell\ll 0, ~\ell\gg 0$, and
$|\ell|\gg 0$. By ${\cal D}^f(R)$ we mean the full subcategory of
${\cal D}(R)$ consisting of complexes $X$ with $\H_{\ell}(X)$ is a
finite $R$--module for all $\ell$.

The left derived functor of the tensor product functor of
$R$-complexes is denoted by $-\otimes_R^{\T{\mathbf L}}-$, and
$\T{\mathbf R}\T{Hom}_R(-,-)$ denotes the right derived functor of
the homomorphism functor of complexes. We need the next two
inequalities for $X,~Y\in{\cal D}_+(R)$ and $Z\in{\cal D}(R)$.

(1.1)\,\,\,\,\,\,\,\, $\inf(X\otimes^{\mathbf L}_RY)\ge \inf
X+\inf Y$\,\,\,\,\, and

(1.2)\,\,\,\,\,\,\, $\sup({\mathbf R}\Hom_R(X,Z))\le\sup Z-\inf
X$.

\noindent For a complex $X$, the dimension of $X$ is defined by
Foxby in [{\bf Fo}] as follows:

(1.3)\,\,\,\,\,\,\, $\dim_RX=\sup\{\dim R/\fp-\inf
X_{\fp}|\fp\in\Spec R\}$.

\noindent It is shown, in [{\bf Fo}; 16.9], that:
$$\dim X = \sup\{\dim_R\H_{\ell}(X)-\ell|\ell\in\Bbb Z \}.$$

\noindent Therefore it is natural to give the following
definition, cf. [{\bf DY}; 2.1]: For a complex $X\in{\cal
D}_+(R)$, the $\fa$--cohomological dimension of $X$ is defined by

$$\cd(\fa,X)=\sup\{\cd(\fa,\H_{\ell}(X))-\ell|\ell\in\Bbb Z\}.$$

\noindent Note that for an $R$--module $M$, this notion agrees
with the classical one.

A complex $X \in {\cal D}_{b} (R) $ is said to be of finite {\it
{projective} }(resp. {\it {injective}}){\it {dimension}} if $ X
\simeq U $, where $U\in{\cal D}_b(R)$ is a complex of projective
(resp. injective) modules.

The full subcategories of $ {\cal D}_{b}(R) $ consisting of
complexes of finite projective (resp. injective) dimension are
denoted by $ {\cal P} (R) $ (resp. $ {\cal I} (R) $). If $X$
belongs to $ {\cal D}_{b} (R) $ , then the following inequalities
hold when $ P \in {\cal P} (R) $ , $ I \in {\cal I} (R) $, cf.
[{\bf Fo}; 8.9 and 8.13].

\[ \begin{array}{l}
\inf ({\mathbf R} \Hom_R (P,X) ) \geqslant  \inf X - \pd _R P;
\\ [.1in]
\inf ({\mathbf R} \Hom_R (X,I) ) \geqslant  - \sup X - \id_R I.
\end{array} \]

\newpage

\noindent{\bf 2. Dimensions and Intersection Theorem}

\vspace{.2in}

In this section we study dimension of the $R$--complex
$\uhom_R(X,Y)$ for appropriate complexes $X$ and $Y$. For
motivation, let us assume $M$ be a finite $R$--module of finite
projective dimension. In [{\bf HZ}] Herzog and Zamani, by using
Buchsbaum--Eisenbud criterion, show that $\dim\Ext^i_R(M,R)\le\dim
R-i$ for all $i\in\Bbb Z$, which is equivalent to say
$\dim\uhom_R(M,R)\le\dim R$. Here, we first bring a generalization
of this result without using of the mentioned criterion; and then
we find some characterizations of complexes in ${\cal D}^f_b(R)$
to have finite projective (or injective) dimensions.

In [{\bf AB}], Auslander and Bridger have introduced the notion of
Gorenstein dimension of a finite $R$--module $M$ as $\gd
M=\sup\{i\in\Bbb Z|\Ext^i_R(M,R)\neq 0\}$, and generalized in
[{\bf Y2}] for reflexive complexes. A complex $X\in{\cal D}(R)$ is
called reflexive if $X\in{\cal}D^f_b(R)$, $\uhom_R(X,R)\in{\cal
D}^f_b(R)$ and the natural map $X\to\uhom_R(\uhom_R(X,R),R)$ is a
quasi-isomorphism. For a reflexive complex $X$, Gorenstein
dimension of $X$ is defined by $\gd X=-\inf\uhom(X,R)$, cf. [{\bf
Y2}; Definition 2.8]. Note that $\gd X\le\pd X$, and equality
holds if $\pd X<\infty$.

\vspace{.1in}

\noindent{\bf Theorem 2.1.} Let $X,Y\in{\cal D}^f_b(R)$. Then
$$\dim\uhom_R(X,Y)\le\dim R+\sup X-\inf Y$$ provided one of the
following conditions holds:

\begin{verse}

(i) $\pd X$ is finite,

(ii) $\id Y$ is finite,

(iii) $\gd X$ and $\pd Y$ are finite.

\end{verse}

\vspace{.1in}

\noindent{\it Proof.} We may assume that $\uhom_R(X,Y)$ is a
homologically non-trivial complex. Choose an integer $\ell$ with
$\H_{\ell}(\uhom_R(X,Y))\neq 0$. Then
$\gr\H_{\ell}(\uhom_R(X,Y))=\depth R_{\fp}$ for some
$\fp\in\Supp\H_{\ell}(\uhom_R(X,Y))$.

(i) Assume $\pd X<\infty$. As
$\H_{\ell}(\uhom_{R_{\fp}}(X_{\fp},Y_{\fp}))\neq 0$, we have
$$-\ell\le -\inf\uhom_{R_{\fp}}(X_{\fp},Y_{\fp})\le\pd
X_{\fp}-\inf Y_{\fp}.$$ Set $\H_{\ell}(\uhom_R(X,Y))=T$. By using
Auslander--Buchsbaum formula, $\pd X_{\fp}=\depth
R_{\fp}-\depth_{R_{\fp}}X_{\fp}$, and so the following computation
hold:

\[ \begin{array}{rl}
\dim T-\ell &\le
\dim T+\pd X_{\fp}-\inf Y_{\fp}\\
&=\dim T+\depth
R_{\fp}-\depth_{R_{\fp}}X_{\fp}-\inf Y_{\fp}\\
&=\dim T+\gr T-\depth_{R_{\fp}}X_{\fp}-\inf
Y_{\fp}\\
&\le\dim R-\depth_{R_{\fp}}X_{\fp}-\inf Y_{\fp}.
\end{array}\]
By [{\bf Fo}; 12.6] we have $-\depth_{R_{\fp}}X_{\fp}\le\sup
X_{\fp}\le \sup X$. Now by using the fact that $\inf Y\le\inf
Y_{\fp}$, we get the assertion.

(ii) Assume $\id Y<\infty$. We have
$$-\ell\le -\inf\uhom_{R_{\fp}}(X_{\fp},Y_{\fp})\le \id
Y_{\fp}+\sup X_{\fp}.$$ Again let $\H_{\ell}(\uhom_R(X,Y))=T$.
Now, using [{\bf Fo}; 13.23 (I')] implies that
\[ \begin{array}{rl}
\dim T-\ell &\le
\dim T+\id Y_{\fp}+\sup X_{\fp}\\
&=\dim T+\depth
R_{\fp}-\inf Y_{\fp}+\sup X_{\fp}\\
&\le\dim T+\gr T-\inf
Y+\sup X\\
&\le\dim R-\inf Y+\sup X.
\end{array}\]
Thus the assertion holds.

(iii) Assume $\gd X<\infty$ and $\pd Y<\infty$. By [{\bf C};
Proposition 2.4.1 and Theorem 2.3.13], we have
\[ \begin{array}{rl}
-\ell\le -\inf\uhom_{R_{\fp}}(X_{\fp}, Y_{\fp})&\le\gd
X_{\fp}-\inf
Y_{\fp}\\
&=\depth R_{\fp}-\depth_{R_{\fp}}X_{\fp}-\inf Y_{\fp}.
\end{array}\]
Now the assertion follows with the same argument as in
(i).\hfill$\square$

\vspace{.2in}

The next Corollary shows that the equality holds for the left
inequality in (T2).

\vspace{.1in}

\noindent{\bf Corollary 2.2} Assume that $M$ is a non--zero finite
$R$--module with $\pd M<\infty$. Then $\dim\uhom_R(M,R)=\dim R$.

\vspace{.1in}

\noindent{\it Proof.} Use Theorem 2.1(i) and (T2).\hfill$\square$

\vspace{.2in}

In [{\bf AB}], Auslander and Bridger show that the ring $R$ is
Gorenstein if and only if all finite $R$--modules have finite
Gorenstein dimension. The following characterization of Gornstein
rings is parallel to the Auslander--Bridger characterization.

\vspace{.1in}

\noindent{\bf Corollary 2.3} Let $(R,\fm,k)$ be a local ring. Then
the following statements are equivalent.

\begin{verse}

(i) $R$ is a Gorenstein ring;

(ii) $\dim\uhom_R(X,R)<\infty$ for all $X\in{\cal D}^f_b(R)$;

(iii) $\dim\uhom_R(M,R)<\infty$ for all finite $R$--module $M$;

(iv) $\dim\uhom_R(k,R)<\infty$.

\end{verse}

\vspace{.1in}

\noindent{\it Proof.} (i)$\Rightarrow$(ii) is clear from Theorem
2.1. (ii)$\Rightarrow$(iii) and (iii)$\Rightarrow$(iv) are
trivial. For (iv)$\Rightarrow$(i), note that
$$\sup\{\ell|\Ext^{\ell}_R(k,R)\neq 0\}=\dim\uhom_R(k,R)<\infty,$$
so $R$ is a Gorenstein ring.\hfill$\square$

\vspace{.2in}

The Bass characterization of Cohen--Macaulay rings says that $R$
is Cohen--Macaulay if and only if there exists a finite
$R$--module with finite injective dimension, cf. [{\bf BH}; Page
375]. The following characterization of Cohen--Macaulay rings is
parallel to the Bass characterization.

\vspace{.1in}

\noindent{\bf Corollary 2.4} Let $(R,\fm,k)$ be a local ring. Then
the following statements are equivalent.

\begin{verse}

(i) $R$ is a Cohen--Macaulay ring;

(ii) There exists a non--zero finite $R$--module $L$ such that
$\dim\uhom_R(X,L)<\infty$ for all $X\in{\cal D}^f_b(R)$;

(iii) There exists a non--zero finite $R$--module $L$ such that
$\dim\uhom_R(M,L)<\infty$ for all finite $R$--module $M$;

(iv) There exists a non--zero finite $R$--module $L$ such that
$\dim\uhom_R(k,L)<\infty$.

\end{verse}

\vspace{.1in}

\noindent{\it Proof.} Note that Cohen--Macaulayness of $R$ is
equivalent to the existence of a non--zero finite $R$--module $L$
of finite injective dimension, and use Theorem 2.1.\hfill$\square$

\vspace{.2in}

Let $(R,\fm)$ be local ring. Let $M$ and $N$ be finite
$R$--modules and $M$ has finite projective dimension. The New
Intersection Theorem of Peskine and Szpiro [{\bf PS}], Hochster
[{\bf H}], and P. Roberts [{\bf R}] yields an inequality
$$\dim N\le\dim (M\otimes_RN)+\pd M.$$ Actually by (T1) we have the
following inequalities,
$$\dim N\le\dim\uhom_R(M,N)\le\dim (M\otimes_RN)+\pd M.$$
On the other hand in [{\bf Fo}] Foxby generalized the New
Intersection Theorem for complexes of modules that is $$\dim
Y\le\dim(X\utp_RY)+\pd X$$ for any $X,Y\in{\cal D}^f_b(R)$ with
$\pd X<\infty$. Now we generalize (T1) for complexes. To do this
we need the following Theorem.

\vspace{.1in}

\noindent{\bf Theorem 2.5}[{\bf Fo}; 9.6 and 11.27(d)] Assume that
$X\in{\cal D}^f_+(R)\cap{\cal P}(R)$ and $Y\in{\cal D}_b(R)$. Then
the following hold.

\begin{verse}

(a) $\uhom_R(X,Y)\simeq\uhom_R(X,R)\utp_RY$

(b) If $R$ is local then $\inf\uhom_R(X,R)=-\pd X$.

\end{verse}

\vspace{.2in}

\noindent{\bf Theorem 2.6} If $X,Y\in{\cal D}(R)$ such that
$\uhom_R(X,Y)$ is not homologically trivial and that $\pd
X<\infty$, then the following statements hold:

\begin{verse}

(i) If $X\in{\cal D}^f_+(R)$ and $Y\in{\cal D}_b(R)$, then
$$\dim\uhom_R(X,Y)\le\dim Y+\pd X.$$

(ii) If $X,Y\in{\cal D}^f_b(R)$, then $$\dim Y+\inf
X\le\dim\uhom_R(X,Y)\le\dim(X\utp_RY)+\pd X+\sup X.$$

\end{verse}

\vspace{.1in}

\noindent{\it Proof.} Choose a prime ideal $\fp\in\Supp X\cap\Supp
Y$ such that

\[\begin{array}{rl}
\dim\uhom_R(X,Y)&=\dim R/\fp-\inf\uhom_R(X,Y)_{\fp}\\
&=\dim R/\fp-\inf\uhom_{R_{\fp}}(X_{\fp},Y_{\fp})\\
&=\dim R/\fp-\inf X^{\star}_{\fp}\utp_{R_{\fp}}Y_{\fp},
\end{array}\]
where $X^{\star}=\uhom_R(X,R)$ and the last equality comes from
Theorem 2.5. Using Nakayama lemma for complexes, we have that

\[ \begin{array}{rl}
\dim R/\fp-\inf\uhom_R(X,Y)_{\fp}&=\dim R/\fp-\inf Y_{\fp}-\inf
X^{\star}_{\fp}\\
&\le\dim Y+\pd X_{\fp}\\
&\le \dim Y+\pd X.
\end{array} \]

(ii). By (i) we have the following inequality
$$\dim\uhom_R(X,Y)\le\dim R/\fp-\inf Y_{\fp}+\pd X$$
for some $\fp\in\Supp X\cap\Supp Y$. We can use Nakayama lemma
again to see that

\[\begin{array}{rl}
\dim\uhom_R(X,Y)&\le\dim R/\fp-\inf(X\utp_R Y)_{\fp}+\inf
X_{\fp}+\pd X\\
&\le\dim(X\utp_RY)+\pd X+\sup X.
\end{array}\]
This proves the right hand side of (ii).

Consider the Intersection Theorem for $X^{\star}$ and $Y$, we have
$$\dim Y\le\dim(X^{\star}\utp_RY)+\pd X^{\star}.$$
Now, by using Theorem 2.5, we get
$$\dim Y\le\dim\uhom_R(X,Y)-\inf X.$$
Now the assertion holds.\hfill$\square$

\vspace{.2in}

\noindent{\bf Corollary 2.7} Let $M$ be a finite $R$--module with
finite projective dimension, and $Y\in{\cal D}^f_b(R)$. Then
$$\dim Y\le\dim\uhom_R(M,Y)\le\dim(M\utp_RY)+\pd M.$$

\newpage

\noindent{\bf 3. Cohomological dimension}

\vspace{.2in}

In [{\bf DY}], the authors investigate the invariants
$-\inf\ugamma_{\fa}(X)$ and $\cd(\fa,X)$, where $\fa$ is an ideal
of $R$, $X\in{\cal D}_+(R)$. The main purpose of this section is
to study $-\inf\ugamma_{\fa}(\uhom_R(X,Y))$ for $X,Y\in{\cal
D}^f_b(R)$. Mainly we seek some results consistent with those in
section 2. First we recall the following Theorem, cf. [{\bf DY};
Theorem 3.2 and Proposition 2.5].

\vspace{.1in}

\noindent{\bf Theorem 3.1} If $X\in{\cal D}_+(R)$ is not
homologically trivial, then
$$-\inf\ugamma_{\fa}(X)\le\cd(\fa,X).$$
Moreover, if $X\in{\cal D}_+^f(R)$, then
\[\begin{array}{rl}
-\inf\ugamma_{\fa}(X)&=\cd(\fa,X)\\
&=\sup\{\cd(\fa,R/\fp)-\inf X_{\fp}|\fp\in\Spec(R)\}.
\end{array}\]

\vspace{.2in}

We are now ready to give a result which is consistent with Theorem
2.6.

\vspace{.1in}

\noindent{\bf Theorem 3.2} Assume $X,Y\in{\cal D}(R)$ and that
$\uhom_R(X,Y)$ is not homologically trivial. The following
statements hold.

\begin{verse}

(i) If $X\in{\cal D}^f_+(R)$ with $\pd X<\infty$ and $Y\in{\cal
D}_b(R)$, then
$$-\inf\ugamma_{\fa}(\uhom_R(X,Y))\le-\inf\ugamma_{\fa}(Y)+\pd X;$$

(ii) If $X,Y\in{\cal D}^f_b(R)$ with $\pd X<\infty$, then
$$-\inf\ugamma_{\fa}(\uhom_R(X,Y))\le-\inf\ugamma_{\fa}(X\utp_RY)+\pd
X+\sup X.$$

\end{verse}

\vspace{.1in}

\noindent{\it Proof.} (i). We have

\[\begin{array}{rl}
-\inf\ugamma_{\fa}(\uhom_R(X,Y))&=-\inf\uhom_R(X,\ugamma_{\fa}(Y))\\
&=-\inf\uhom_R(X,R\utp_R\ugamma_{\fa}(Y))\\
&=-\inf\uhom_R(X,R)\utp_R\ugamma_{\fa}(Y)\\
&\le-\inf X^{\star}-\inf\ugamma_{\fa}(Y),
\end{array} \]
where the first equality holds by [{\bf Y1}; Proposition 2.3] and
the third one is Tensor-Hom evaluation. Now the proof of (i) is
finished.

(ii). Assume $X,Y\in{\cal D}^f_b(R)$, so that
$\uhom_R(X,Y)\in{\cal D}^f_+(R)$ and thus, by Theorem 3.1 and the
proof of Theorem 2.6(i),
$$-\inf\ugamma_{\fa}(\uhom_R(X,Y))=\cd(\fa,R/\frak
p)-\inf X^{\star}_{\fp}-\inf Y_{\fp}$$ for some $\fp\in\Supp
X\cap\Supp Y$. Now with the same argument as in Theorem 2.6(ii)
and using Theorem 3.1 again, we get the result.\hfill$\square$

\vspace{.2in}

To prove a similar statement as left side equality of Theorem
2.6(ii), one need to show the Intersection Theorem for local
cohomology. That is if $X,Y\in{\cal D}^f_b(R)$ and $\pd X<\infty$,
then
$$-\inf\ugamma_{\fa}(Y)\le-\inf\ugamma_{\fa}(X\utp_RY)+\pd X.$$ Although the proof is not known in
general, for the authors, it may worth noting that there is a
proof in case $\Supp Y\subseteq\Supp X$. To proceed, it is easy to
show that if $X\in{\cal D}^f_b(R)$ with finite projective
dimension, then $\Supp X^{\star}\subseteq\Supp X$. Thus $\Supp
X=\Supp X^{{\star}{\star}}\subseteq\Supp X^{\star}\subseteq\Supp
X$ gives $\Supp X^{\star}=\Supp X$. Now, we can state the
following result.

\vspace{.1in}

\noindent{\bf Theorem 3.3} Assume $X,Y\in{\cal D}^f_b(R)$ with
$\pd X<\infty$ and $\Supp Y\subseteq\Supp X$. Then

\begin{verse}

(i) $-\inf\ugamma_{\fa}(Y)\le -\inf\ugamma_{\fa}(X\utp_RY)+\pd X.$

(ii) $-\inf\ugamma_{\fa}(Y)+\inf X\le
-\inf\ugamma_{\fa}(\uhom_R(X,Y)).$

\end{verse}

\noindent{\it Proof.} (i). By using Theorem 3.1, there exists
$\fp\in\Supp Y$ such that
$$-\inf\ugamma_{\fa}(Y)=\cd(\fa,R/\fp)-\inf Y_{\fp}.$$
As $\Supp Y\subseteq\Supp X$, then

\[\begin{array}{rl}
-\inf\ugamma_{\fa}(Y)&=\cd(\fa,R/\fp)-\inf Y_{\fp}\\
&=\cd(\fa,R/\fp)-\inf(X\utp_RY)_{\fp}+\inf X_{\fp}\\
&\le -\inf\ugamma_{\fa}(X\utp_RY)+\pd X.
\end{array}\]
The proof of (i) is finished.

(ii) As $\Supp Y\subseteq\Supp X^{\star}$, we have

\[\begin{array}{rl}
-\inf\ugamma_{\fa}(Y)&\le -\inf\ugamma_{\fa}(X^{\star}\utp_RY)+\pd X^{\star}\\
&=-\inf\ugamma_{\fa}(\uhom_R(X,Y))-\inf X,

\end{array}\]
where the equality holds by Theorem 2.5(a) and so the assertion
holds.\hfill$\square$

\vspace{.2in}

Now, we are able to bring the following result.

\vspace{.1in}

\noindent{\bf Corollary 3.4} Assume that $M$ is a finite
$R$--module with finite projective dimension. Then, for any
$Y\in{\cal D}^f_b(R)$ with $\Supp Y\subseteq\Supp M$,
$$-\inf\ugamma_{\fa}(Y)\le-\inf\ugamma_{\fa}(\uhom_R(M,Y)).$$

In particular, if $\Supp M=\Spec R$, then
$$-\inf\ugamma_{\fa}(R)\le-\inf\ugamma_{\fa}(\uhom_R(M,R)).$$

\vspace{.3in}

\noindent {\large\bf Acknowledgments.} The authors would like to
thank the referee for his/her substantial comments.

\baselineskip=16pt

\begin{center}
\large {\bf References}
\end{center}
\vspace{.2in}

\begin{verse}

[A] M. Auslander, {\em Anneaux de Gorenstein, et torsion en
alg\`{e}bre commutative}, (French) S\'{e}minaire d'Alg\`{e}bre
Commutative dirig\'{e} par Pierre Samuel, 1966/67. Texte
r\'{e}dig\'{e}, d'apr\`{e}s des expos\'{e}s de Maurice Auslander,
Marquerite Mangeney, Christian Peskine et Lucien Szpiro. \'{E}cole
Normale Sup\'{e}rieure de Jeunes Filles Secr\'{e}tariat
math\'{e}matique, Paris 1967.

[AB] M. Auslander; M. Bridger, {\em Stable module theory}, Memoirs
of the American Mathematical Society, No. 94 American Mathematical
Society, Providence, R.I. 1969.

[BH] W. Bruns; J. Herzog, {\em Cohen-Macaulay rings}, Cambridge
Studies in Advanced Mathematics, {\bf 39}. Cambridge University
Press, Cambridge, 1993.

[C] L. W. Christensen, {\em Gorenstein dimensions}, Lecture Notes
in Mathematics, {\bf 1747} Springer-Verlag, Berlin, 2000.

[DY] M. T. Dibaei; S. Yassemi, {\em Cohomological dimension of
complexes}, to appear in Comm. Algebra.

[Fa] G. Faltings, {\em \"{U}ber lokale Kohomologiegruppen hoher
Ordnung}, J. Reine Angew. Math. {\bf 313} (1980), 43--51.

[Fo] H.-B. Foxby, {\em  Hyperhomological algebra and commutative
algebra}, notes in preparation.

[H] M. Hochster, {\em Topics in the homological theory of modules
over commutative rings}, Regional Conference Series in
Mathematics, No. 24. American Mathematical Society 1975.

[HL] C. Huneke; G. Lyubeznik, {\em On the vanishing of local
cohomology modules}, Invent. Math. {\bf 102} (1990), 73--93.

[HZ] J. Herzog; N. Zamani, {\em Duality and vanishing of
generalized local cohomology}, preprint 2003.

[PS] C. Peskine; L. Szpiro, {\em Dimension projective finie et
cohomologie locale}, Inst. Hautes Études Sci. Publ. Math. No. 42
(1973), 47--119.

[R]  P. Roberts {\em Le th\'{e}or\`{e}me d'intersection}, C. R.
Acad. Sci. Paris S\'{e}r. I Math. 304 (1987), 177--180.

[Y1] S. Yassemi, {\em Generalized section functors}, J. Pure Appl.
Algebra {\bf 95} (1994), 103--119.

[Y2] S. Yassemi, {\em G--dimension}, Math. Scand. {\bf 77} (1995),
161--174.

\end{verse}

\end{document}